\documentclass{amsart}
\usepackage{graphicx,amsmath,amssymb,amsfonts,amsthm}
\usepackage[all,cmtip]{xy}
\usepackage{hyperref}
\newtheorem{thm}{Theorem}[section]

\theoremstyle{definition}

\theoremstyle{remark}
\newtheorem{rem}[thm]{Remark}
\numberwithin{equation}{section}

\newcommand{\Spec}{ \mathrm{Spec}}

\newcommand{\Mod}{\mathsf{Mod}}

\newcommand{\Pic}{\mathrm{Pic}}

\newcommand{\Fr}{F}
\newcommand{\Xp}{X^{(p)}}

\newcommand{\OO}{\mathcal{O}}
\newcommand{\PP}{\mathbb{P}}

\newcommand{\ZZ}{\mathbb{Z}}
\newcommand{\QQ}{\mathbb{Q}}

\newcommand{\ZZpurcom}{\widehat{\mathbb Z_p^{\mathrm{ur}}}}
\newcommand{\FFpbar}{\overline{\mathbb{F}}_p}
\newcommand{\Proj}{\mathrm{Proj}}

\newcommand{\phirel}{\phi_{X/S}}
\newcommand{\Frel}{F_{X_0/S_0}}
\newcommand{\phiabs}{\phi_{X}}

\newcommand{\Crel}{C_{X_0/S_0}}

\newcommand{\Xphi}{X^{(\phi)}}
\newcommand{\Hunder}{\underline{H}}
\newcommand{\id}{\mathrm{id}}
\newcommand{\rk}{\mathrm{rk}}

\title{Positivity and Lifts of the Frobenius}
\author{Taylor Dupuy}
\address{Department of Mathematics and Statistics, MSC01 1115, 1 University of New Mexico, Albuquerque, New Mexico, 87131-0001}
\email{dupuy@math.unm.edu}

\begin{document}

\begin{abstract}
If $X/W(\FFpbar)$ is a smooth projective scheme with positive Kodaira dimension then $X$ does not have a lift of the $p$-Frobenius.
\end{abstract}

\maketitle

\section{Introduction }\label{sec:introduction}
This paper is motivated by J. Borger's work on $\Lambda$-schemes \cite{Borger2009}. A $\Lambda$-scheme is a scheme $X/\ZZ$ with a commuting family of lifts of the absolute $p$-Frobenius for each prime $p$ and according to Borger the lifts should be viewed as descent data a category of ``schemes over the field with one element'' (think Weil-Descent). 

Throughout this paper $p$ will be a fixed prime and $W=W(\FFpbar)$ will be the ring of $p$-typical Witt vectors of $\FFpbar$. 
Recall that it is the unique $p$-adically complete discrete valuation ring that it is unramified over $\ZZ_p$ and has residue field $\FFpbar$. 
In this paper we show that flat projective schemes $X/W$ of positive Kodaira dimension cannot have a lift of the Frobenius. 
We do this by generalizing a result of Raynaud to higher dimensions. 

Let's fix the following notation: If $A$ is a flat $W$-algebra then we write
 $$A_n:= A/p^{n+1}A$$
If $X/W$ is a flat scheme then we use the similar notation
  $$ X_n = X \otimes_W W_n.$$ 
If $A$ is a $W$-algebra, a \textbf{lift of the absolute Frobenius} on $A/pA$ is a ring endomorphisms $\phi: A\to A$ with the property that for all $a\in A$ we have
 $$ \phi(a) \equiv a^p \mod pA.$$
The ring $W$ has a unique lift of the Frobenius $\phi$. If we present $W$ as 
\begin{equation}
W=\ZZ_p[\zeta : \zeta^n = 1, p \nmid n]^{\widehat \ }, 
\end{equation}
then $\phi(\zeta) = \zeta^p$ for all roots of unity and $\phi\vert_{\ZZ_p} = \id$.
The hat in the above displayed equation denotes $p$-adic completion.

Every scheme $X_0/\FFpbar$ admits sheaf endomorphism $F^*:\OO_{X_0} \to \OO_{X_0}$ called the \textbf{absolute Frobenius} such that $F^*(s) = s^p$ for a local section $s\in \OO_{X_0}$. This sheaf endomorphism defines a unique morphisms of schemes $F:X_0\to X_0$ which is constant on the underlying topological space and acts as stated on local sections of the structure sheaf. If $X/W$ is a flat scheme which reduces to $X_0$ mod $p$. 
We call a morphism of sheaves $\phi_X: X\to X$ a \textbf{lift of the Absolute Frobenius} provided it is equal to the absolute Frobenius mod $p$. 
We discuss lifts of the Frobenius further in section \ref{sec:frobenius}.

This paper generalizes what is known in the case of (smooth projective) curves which can be described by the table below:

\begin{equation}
 \begin{tabular}{|c||c|c|c|}
\hline Genus & $g(C)=0$ & $g(C)=1 $ & $g(C)\geq 2$ \\
\hline Has lift? & Always & Rarely & Never \\
\hline 
\end{tabular}  
\label{eqn:curves}
\end{equation}

For the first column $C = \PP^1_W$ and we can construct the lift by taking $p$-th powers on the coordinates and patching the map together.  
By this we mean that we cover $\PP^1$ with two copies of $\AA^1 = \Spec(W[T])$ in the standard way and define the lift $\phi_{\PP^1}$ on these affine open sets by $\phi_{\PP^1}^*: T \mapsto T^p$. 
One can check that this gives a well-defined lift of the Frobenius on the intersection of the two open sets. 
The elliptic curve case is essentially the theory of the canonical lifts due to Serre and Tate and is reviewed in Katz's classic article \cite{Katz1981}. See section 4.1 of \cite{Buium2009} for a summary.
The hyperbolic case is due to Raynaud and is a Cartier operator argument. Naively, one can reformute the \ref{eqn:curves} in terms of Kodaira dimension (see section \ref{sec:kodaira dimension}) and ask if a similar table holds. In the case of curves $g=0$ if and only if $\kappa(C)<0$, $g=1$ if and only if $\kappa(C)=0$ and $g\geq 2$ if and only if $\kappa(C)=1$. In higher dimensions one finds:
\begin{equation}\label{eqn:varieties}
 \begin{tabular}{|c||c|c|c|}
\hline Kodaira dimension & $\kappa(X)<0$ & $\kappa(X)=0$ & $\kappa(X)>0$   \\
\hline Has lift? & Very rarely & Rarely & Never \\
\hline 
\end{tabular}  
\end{equation} 

This main result of this paper is a modest modernization of Raynaud's Lemma which appears to be missing in the literature.
\begin{thm}[Higher Dimensional Raynaud's Lemma]\label{thm:main}
Let $W=\ZZpurcom$ and $X/W$ be a smooth projective scheme. If $X$ has positive Kodaira dimension then $X$ does not admit a lift of the Frobenius.
\end{thm}
The strategy of the above version of Raynaud's Lemma (\cite{Raynaud1983} Lemma I.5.1) is roughly the same as the original: Supposing a lift of the Frobenius, cook up a map between the de Rham complexes and show that it is non-zero using the theory of the Cartier operator. 
Next, using the positivity of the Kodaira Dimension, show that such a non-zero map can not exist giving a contradiction.

 It is natural to ask whether or not every Fano variety admits a lift of the Frobenius. 
The answer to this question is negative. 
Counter examples can be found via a paper of Paranjape and Srinivas \cite{Paranjape1989}.
In particular the remarks following problem 1 on page 1 of \cite{Paranjape1989} state the following: 
If $G/K$ a semisimple algebraic group over an algebraically closed field $K$ of characterisic zero and $P$ is a parabolic subgroup, then the variety $Y := G/P$ does not admit a self map $f:Y \to Y$ of degree bigger than one unless $Y \cong \PP^n_K$. 
Since Grassmannians take this form and all Grassmannians are Fano and defined over (say) $W(\FFpbar)$ there do not exists lifts of the Frobenius $Y$ where the generic fiber is not isomorphic to $\PP^n_K$. 
Otherwise one could obtain a contradiction after base change.
This was found in section 2.8 of Borger's paper \cite{Borger2009}. 
Here they show that the only flag variety which admits a lift of the Frobenius is projective space.

In section \ref{sec:background} we provide the necessary background and in section \ref{sec:proof} we give the proof.

Conversations with Alexandru Buium, James Borger and Michael Nakamaye were very helpful in preparing this paper.

\section{Background Necessary for the Proof}\label{sec:background}

Section \ref{sec:frobenius} reviews the notions of absolute and relative Frobenius morphisms, section \ref{sec:cartier operator} recalls the Cartier operator, and section \ref{sec:kodaira dimension} recalls the notion of Kodaira dimension.

\subsection{Lifts of the Frobenius}\label{sec:frobenius}
Let $S_0 = \Spec(\FFpbar)$ and $F_{S_0}$ denote the absolute Frobenius on $S_0$. If $u_0:X_0\to S_0$ is a morphism of schemes we have the following commutative diagram which defines the scheme $X_0^{(p)}$ and the maps $\Fr_{X_0/S_0}$ and $\varphi_0$:
	\begin{equation}\label{eqn:frobenius_diagram}
	\xymatrix{
	X_0\ar@{-->}[dr]^{\Fr_{X_0/S_0}} \ar@/^2pc/[drr]^{\Fr_{X_0}} \ar@/_2pc/[ddr]_{u_0} & &\\ 
	& X_0^{(p)}  \ar[r]^{ \ \ \ \ \ \ \ \ \ \ \ \varphi_0} \ar[d]^{u_0^{(p)}} & X_0\ar[d]^{u_0}\\
	& S_0 \ar[r]^{\Fr_{S_0} } & S_0\\
	}
	\end{equation}
In the above diagram we have $X_0\times_{S_0,\Fr_{S_0}}S_0$. The morphism $\Fr_{X_0/S_0}$ is called a \textbf{relative Frobenius}. For all $f$ local sections of $\OO_{X_0}$ and all $c\in \FFpbar$ we have
\begin{eqnarray*}
 \varphi_0^*(f) &=& f \otimes 1, \\
 \varphi_0^*(cf) &=& c^p \varphi_0^*(f).
\end{eqnarray*}
Commutivity of the diagram implies that 
\begin{equation}
 F_{X_0/S_0}^*(\varphi^*_0 f) = \Fr_{X_0/S_0}^*(f\otimes 1) = \Fr_{X_0}^*(f) = f^p.
\end{equation}
Commutivity also says that unlike $F_{X_0}$ the map $\Fr_{X_0/S_0}^*:\OO_{\Xp_0} \to \OO_{X_0}$ defines a morphism of schemes over $\FFpbar$. For all $c\in \FFpbar$ and $f$ local sections of $\OO_{X_0}$ 
\begin{eqnarray*}
 c \cdot \varphi_0^*( f ) &=& f \otimes c,\\ 
F_{X_0/S_0}^*(c \cdot \varphi^*_0(f) ) &=& cf^p .
\end{eqnarray*}
We employ the notation
\begin{equation}
f^{(p)} := \varphi_0^*(f) = f\otimes 1 ,
\end{equation}
so that for all $f\in \OO_{X_0}$ and all $c\in \FFpbar$
\begin{eqnarray*}
(cf)^{(p)} &=& c^p f^{(p)}, \\
F_{X_0/S_0}^* (f^{(p)}) &=& f^p.
\end{eqnarray*}

Suppose now that $X/S$ is flat with reduction mod $p$ equal to $X_0/S_0$ then a \textbf{lift of the absolute Frobenius} is a map $\phi_{X}: X\to X$ whose reduction mod $p$ induces $\Fr_{X_0}$. 

Suppose $\phi: S\to S$ is a lift of the absolute Frobenius.
Let $X^{\phi} = X \times_{S,\phi} S$ be the pullback of $X$ by $\phi$ with $\varphi: X^{(\phi)} \to X$ given by the first projection. 
We call $\phi_X$ a \textbf{lift of the absolute Frobenius} (compatible with the lift $\phi$ on $S$) and $\phi_{X/S}$ a \textbf{ lift of the relative Frobenius} on $X$ (compatible with the lift $\phi$ on $S$) if we have the following diagram:
\begin{equation}\label{eqn:frobenius_lift_diagram}
	\xymatrix{
	X\ar[dr]^{\phi_{X/ S}} \ar@/^1pc/[drr]^{\phi_{X}} \ar@/_1pc/[ddr]_{u} & &\\ 
	& X^{(\phi)} \ar[r]^{\varphi} \ar[d]^{u^{(\phi)}} & X\ar[d]^{ u}\\
	& S \ar[r]^{\phi} & S\\
	}.
\end{equation}
Note from the above diagram that lift of the absolute Frobenius on $X$ exists if and only if a lift of the relative Frobenius exists.

\subsection{Cartier Operator}\label{sec:cartier operator}

Let $F$ be a sheaf on a scheme $X$. In what follows we let $\Hunder^{\bullet}(F)$ denote the sheaf associated to the presheaf $U \mapsto H^{\bullet}(U,F)$. 
Also we will let ``$ s\in F $''  denote``$s\in F(U)$ for some open subset $U$ of $X$''.
 
\begin{thm}[Cartier]
 Let $X_0/S_0$ be a smooth projective scheme where $S_0/\FFpbar$. There exists a unique isomorphism of graded $\OO_{X^{(p)}_0}$-algebras 
 	$$C_{X_0/S_0}^{-1}: \bigoplus_{i} \Omega^i_{X^{(p)}_0/S_0} \to \bigoplus_{i} \Hunder^i( F_{X_0/S_0*}\Omega_{X_0/S_0}^{\bullet})$$
 such that for every $\xi ,\eta \in \bigoplus_{i} \Omega^i_{X^{(p)}_0/S_0}$ and $f\in \mathcal O_{X^{(p)}}$ we have 
 \begin{enumerate}
  \item $C^{-1}_{X_0/S_0}(1) = 1$
  \item $C^{-1}_{X_0/S_0}(\xi \wedge \eta) = C^{-1}_{X_0/S_0}(\xi) \wedge C_{X_0/S_0}^{-1}(\eta) $
  \item $C^{-1}_{X_0/S_0}( d\varphi^* f ) = [f^{p-1}df]$ in $\Hunder^1(F_{X_0/\FFpbar *} \Omega^{\bullet}_{X_0/\FFpbar})$
 \end{enumerate}
 and $C_{X_0/S_0}^{-1}$ is an isomorphism.
\end{thm}
For a proof of this theorem see \cite{Katz1970} Theorem (7.2).

\begin{rem}
The map $C^{-1}_{X_0/S_0}$ in the above proposition is called the \textbf{inverse of the Cartier operator}. The \textbf{Cartier operator} $C_{X_0/S_0}: \bigoplus_{i} \Hunder^i( F_{X_0/S_0*}\Omega_{X/S}^{\bullet}) \to \bigoplus_{i} \Omega^i_{X^{(p)}_0/S_0}$ is the inverse of $C^{-1}_{X_0/S_0}$. On each component they are both morphisms of sheaves of $\OO_{X_0^{(p)}}$-modules.
\end{rem}

\subsection{Kodaira Dimension}\label{sec:kodaira dimension}

Let $X/W$ be a smooth scheme. Then $X$ admits a canonical sheaf $\omega$ and one can define the \emph{canonical ring} as
\begin{equation}
R(X,\omega) := \bigoplus_{n\geq 0} H^0(X,\omega^{\otimes n})
\end{equation}
which is a graded ring which has the structure of $W$-algebra.

One should recall that for every $n$ the module $H^0(X,\omega^{\otimes n})$ is free and finite rank since $W$ is local.

We should recall that for every $X$ there exists some polynomial $P_X(t)$ which takes integral values and some integer $N\geq 0$ such that for all $n>N$ 
 $$h^0(X,\omega^{\otimes n}) = P_X(n).$$ 
Here $h^0(X,\omega^{\otimes n}) = \mathrm{rk}_W H^0(X,\omega^{\otimes n})$. Note that this definition makes sense since $H^0(X, \omega^{\otimes n})$ is free:
Since $H^0(X,\omega^{\otimes n}) = f_*\omega^{\otimes n}_{X/S}(X)$ where $f: X \to S = \Spec(W)$ is the structure sheaf and $f_*\omega^{\otimes n}_{X/S}$ is locally free on $S$, this $H^0(X,\omega^{\otimes n})$ a locally free $\OO_X$-module. But since $W$ is local flatness, freeness and local freeness are all the same thing and $H^0(X,\omega^{\otimes n})$ is free of finite rank. 
The statement about the Hilbert polynomial actually being a polynomial then follows from Hilbert's syzygy theorem over fields, Nakayama's Lemma and flatness.

The \textbf{Kodaira dimension} (introduced by Iitaka) of $X$ is defined to be the dimension of the projective scheme associated to canonical ring:
 \begin{equation}
  \kappa(X) = \dim_W \Proj \ R(X,\omega).
 \end{equation}
Here $\dim_W$ denotes the relative dimension of the scheme over $W$. 
There are many equivalent ways to define the Kodaira Dimension. 
We will make use of the equivalent definition 
  $$ \kappa(X) = \inf \lbrace  a : \exists N\geq 0, \exists C\geq 0, \forall d \geq N, \frac{\rk_W (H^0(X,\omega^{\otimes d}). }{d^a} < C \rbrace. $$

\begin{rem}
 We will find it convenient to use divisorial notation in place of sheaf notation in places. 
 
 \begin{itemize}
  \item In this case $K$ will denote a (fixed representative of) the canonical divisor and $\omega \cong \OO(K)$.
  \item One should also recall that $\OO(nD) \cong \OO(D)^{\otimes m}$ for integers $m$ where $L^{\otimes -1}$ will denote the dual of a sheaf of modules $L$.
 \end{itemize}

\end{rem}

 Let $X/L$ be a smooth schemes over a field $L$ of characteristic zero. The condition $\kappa(X)>0$ is equivalent to the canonical divisor $K$ being rationally equivalent to a non-trivial divisor with non-negative $\QQ$-coefficients

 Here is how to see this: by second definition of Kodaira dimension and the assumption that $\kappa(X) > 0$ implies that there exists some $m$ and some non-zero $s \in H^0(X,\OO(mK))$. By definition of $\OO(mK)$ we have $ D:= (s) + mK \geq 0$ which proves that $\frac{1}{m}D \sim K$ as $\QQ$-divisors.
 
 Conversely, suppose that $\frac{1}{m}D \sim K$. This means that $D \sim mK$ and $H^0(X_L,\OO(mK)) = H^0(X, \OO(D))$ where $D$ is an effective divisor with $\ZZ$ coefficients. Since $D$ is effective it has a nonconstant section $s\in \OO(D)(X)$ which tells us that $L \oplus L s  \subset \OO(D)(X).$ By looking at powers of $s$ this we have that the dimension of $\OO(nD)(X)$ is bigger than $n+1$. 
 Hence we have $\dim_L \OO(nmK)(X)/(nm)^a \geq (n+1)/(nm)^a$ and in order for this to be bounded as $n\to \infty$ we need $a\geq 1$. 

 \section{Proof of Main Theorem}\label{sec:proof}

The aim of this section is to prove theorem \ref{thm:main}. Here we suppose that $X/S$ is a smooth projective scheme where $S = \Spec(W)$ which admits lifts $\phirel:X \to \Xphi$ and $\phi_X$ of the relative and absolute Frobeniuses on $X_0$. Note that these are automatically compatible with the lift of the Frobenius on $S$, $\phi_S$, since it is unique.

Our plan is use the Frobenius lift $\phi_X$ to construct a morphism $\tau^{\bullet}$ which lifts the inverse of the Cartier operator and derive a contradiction. 
We will often write out our formulas using the absolute Frobenius and the corresponding arrows with the relative Frobenius.  

For local sections $x$ of $\OO_X$ we can write
 \begin{equation}
  \phirel^*( x\otimes 1)= \phiabs^*(x)  = x^p + p \delta(x)
 \end{equation}
for some $\delta(x) \in \OO_X$. Note that the above expression defines the morphisms of sheaves of sets $\delta: \OO_{\Xphi} \to \phi_{X/S*}\OO_X$. 
Since $X$ is flat over $W$, multiplication by $p$ is injective and division by $p$ makes sense so that we have $\delta(x) = (\phiabs^*(x) - x^p)/p$ for a local section of $s\in \OO_X$. 

On differentials $dx \in \Omega_{X/S}$ or $d(x\otimes 1) \in \Omega_{\Xphi/S}$ we have 
  \begin{equation}
  \phirel^*( d x \otimes 1 )= \phiabs^* dx =  p x^{p-1} dx + p d \delta(x) \in \phi_{X/S*} \Omega_{X/S}.
  \end{equation}
This implies that elements in the image of the map $\phi_{X/S}^*: \Omega_{X^{(\phi)}/S} \to \phi_{X/S*}\Omega_{X/S}$ (equivalently, the map $\Omega_{X/S} \to \phi_{X*}\Omega_{X/S}$) are divisible by $p$. We extend this map multiplicatively: For $\omega = dx_1 \wedge \cdots \wedge dx_r$ a local section of the $r$th wedge power of the module of Kahler differentials $\bigwedge^r \Omega_{X/S}$ we have 
 \begin{equation}
  \phiabs^*( dx_1 \wedge \cdots \wedge dx_r ) \in p^r \Omega_{X/S}
 \end{equation}
at each stage if we divide by the appropriate power of $p$ ($p^r$ for $r$-forms) we get a morphism of complexes $\tau^{\bullet}: \Omega_{\Xphi/S}^{\bullet} \to (\phirel)_* \Omega_{X/S}^{\bullet}$. More precisely $\tau^i: \Omega^i_{\Xphi/S} \to \phi_{X/S*}\Omega^i_{X/S}$ is defined by $\tau^i = \frac{1}{p^i}\phi_{X/S}^*$ for $i = 1, \ldots, \dim_S(X)$. Note that $\tau^1(dx) = x^{p-1}dx + d\delta(x)$ which a lift of the inverse of the Cartier operator (in the sense that it is a morphism of modules that reduces to the $C^{-1}_{X_0/S_0}$). 

Observe that the image of $\tau^{\bullet}$ is closed in the sense that for each $i$ and every local section $s \in \Omega^i_{\Xphi/S}$ we have $d(\tau^i(s)) =0$. 
If we let $\pi^{\bullet}$ denote the map which takes closed forms to their class in cohomology we have 
\begin{equation}\label{eqn:get cartier}
 \Crel^{\bullet} \circ (\pi^{\bullet} \circ \tau^{\bullet})_0 = \mbox{identity},
\end{equation}
where $\Crel^{\bullet}$ is the Cartier isomorphism. Equation \ref{eqn:get cartier} just implies that $(\pi^{\bullet}\circ \tau^{\bullet})_0$ is just the inverse of the Cartier isomorphism defined in section \ref{sec:cartier operator}. This implies that $(\pi^{\bullet} \circ \tau^{\bullet})_0$ is non-zero which implies that $\pi^{\bullet}\circ \tau^{\bullet}$ is non-zero which implies that $\tau^{\bullet}$ is non-zero. Define $(\tau')^i$ to be the adjoint of $\tau^i$ for each $i$. By adjointness we also have that $(\tau')^{\bullet}: \phirel^{*} \Omega_{\Xphi/S}^{\bullet} \to \Omega_{X/S}^{\bullet} $ is nonzero; 
in particular since $\tau^i$ is nonzero for each $0\leq i\leq \dim_S(X)$ we have $\tau^{i}{}'$ being non-zero for each $0 \leq i \leq \dim_S(X)$. 
More generally, since $(\pi^{\bullet} \circ \tau^{\bullet})_0 = \pi^{\bullet}_0 \circ \tau_0^{\bullet}$ we have that $\tau_0'{}^{\bullet}$ is nonzero in each degree. 

We will now derive our contradiction by examining the the top degree piece: In degree $m = \dim(X)$ we have a non-zero morphisms of $\OO_{X_0}$-modules
  \begin{equation}
   (\tau_0')^m: \Frel^*\omega_{X_0^{(p)}/S_0} \to \omega_{X_0/S_0}.
  \end{equation}
Here we just used that the canonical sheaf $\omega_{X_0/S_0}$ is equal to $\Omega^m_{X_0/S_0}$ since $X/S$ is smooth. 
We claim that such a non-zero map can't exist. 
In terms of divisors these line bundles are
\begin{eqnarray*}
 \Frel^*\omega_{X_0^{(p)}/S_0} \cong \omega_{X_0/S_0}^{\otimes p} \cong \OO(pK_{X_0}), & \omega_{X_0/S_0} \cong \OO (K_{X_0} ).
\end{eqnarray*}
The isomorphisms $F_{X_0/S_0}^* \omega_{X_0^{(p)}/S_0} \cong \omega_{X_0/S_0}^{\otimes p}$ follows from considered the affect of pulling back a line bundle by the Frobenius on the transition data in $H^1(X_0, \OO_{X_0}^{\times}) = \Pic(X_0)$. 
Since $(\tau_0')^m \neq 0$ we have 
\begin{eqnarray*}
 0 \neq (\tau_0')^m \in \Mod_{X_0}(\OO(pK_{X_0}),\OO(K_{X_0})) &\cong& \Mod_{X_0}(\OO,\OO( (1-p)K_{X_0}))\\
 &\cong& H^0(X_0,(1-p)K_{X_0})
\end{eqnarray*}
is nonzero (here $\Mod_X$ denotes the category of $\OO_X$-modules and $\Mod_X(A,B)$ denotes the morphisms of $A\to B$ of $\OO_X$-modules). By flatness of $X/W$ we have $h^0(X_0,(1-p)K_{X_0})=h^0(X_L,(1-p)K_{X_L})$ where $L = \mathrm{Frac}(W) = \widehat{\QQ}^{\mathrm{ur}}_p$ and $X_L$ is the generic fiber of $X/W$. Since $\kappa(X_L)>0$ is equivalent to $K_{X_L}\sim_{\QQ} D \geq 0$ as a $\QQ$-divisor we have $h^0(X_L,(1-p)K_{X_L})=0$ which gives a contradiction.

\bibliographystyle{alpha}
\bibliography{../../bib/clean}

\end{document}